 \documentclass[preprint,1p,times]{elsarticle}

\usepackage{amssymb}
\usepackage{color} % statements \color
\def\Universite{Universite\ }
\def\Mathematiques{Mathematiques\ }
\def\Departement{Departement\ }
\def\Rivieres{Rivieres\ }
\def\Montreal{Montreal\ }
\def\Quebec{Quebec\ }
\def\agrave{a\ }

\def \ccomma{\raise 2pt\hbox{,}\ } % Le petit livre de TeX page 234
\def \D {\hbox{d}}

\def \metric{\upsilon}

\def \mod#1{\vert #1 \vert}
\def \barQ {\overline{Q}}
\def \barz {{\bar z}}
\def \cz {a}

\def \modQ {\mod{Q}}

\def \PVI    {{\rm P_{\rm VI}}}
\def \PV     {{\rm P_{\rm V}}}
\def\Fa{F_a}
\def\Fb{F_b}

\def\barz{{\bar z}}

\journal{Journal of geometry and physics}

\begin{document}

\begin{frontmatter}

%% Title, authors and addresses

%% use the tnoteref command within \title for footnotes;
%% use the tnotetext command for theassociated footnote;
%% use the fnref command within \author or \address for footnotes;
%% use the fntext command for theassociated footnote;
%% use the corref command within \author for corresponding author footnotes;
%% use the cortext command for theassociated footnote;
%% use the ead command for the email address,
%% and the form \ead[url] for the home page:
%% \title{Title\tnoteref{label1}}
%% \tnotetext[label1]{}
%% \author{Name\corref{cor1}\fnref{label2}}
%% \ead{email address}
%% \ead[url]{home page}
%% \fntext[label2]{}
%% \cortext[cor1]{}
%% \affiliation{organization={},
%%             addressline={},
%%             city={},
%%             postcode={},
%%             state={},
%%             country={}}
%% \fntext[label3]{}
%
%
\title{A purely analytic derivation of Bonnet surfaces}
%\title{The eleven meromorphic traveling waves of cubic and quintic complex Ginzburg-Landau equations}
%
%use optional labels to link authors explicitly to addresses:
% ERROR of elsarticle.cls if blanks like in \author[n ]
\author[Borelli,HKUMath]{Robert Conte\corref{cor1}}\ead{Robert.Conte@cea.fr,         ORCID https://orcid.org/0000-0002-1840-5095}
\author[CRM,UQTR]       {Alfred Michel Grundland}  \ead{Grundlan@CRM.UMontreal.CA,   ORCID https://orcid.org/0000-0003-4457-7656}
\cortext[cor1]{Corresponding author.
Robert Conte,
Centre Borelli,
\'Ecole normale sup\'erieure Paris-Saclay,
4 avenue des sciences,
F-91190 Gif-sur-Yvette,
France.
Robert.Conte@cea.fr
}
%
%affiliation[Borelli]{organization={Universit\'e Paris-Saclay, ENS Paris-Saclay, CNRS}, \'e => 6 errors of elsarticle.cls
\affiliation[Borelli]{organization={\Universite Paris-Saclay, ENS Paris-Saclay, CNRS},
            addressline={Centre Borelli},
            city={Gif-sur-Yvette},
            postcode={91140},
           %state={},
            country={France}}
\affiliation[HKUMath]{organization={The University of Hong Kong},
            addressline={Department of Mathematics, Pokfulam},
           %city={Hong Kong},
           %postcode={},
           %state={},
            country={Hong Kong}}
\affiliation[CRM]{organization={Centre de Recherches \Mathematiques},
            addressline={\Universite de \Montreal, C.P. 6128, Succ. Centre-ville},
            city={\Montreal}, % 
            postcode={H3C 3J7},
           %state={},
            country={\Quebec, Canada}}
\affiliation[UQTR]{organization={\Universite du \Quebec \agrave Trois-{\Rivieres}},
            addressline={\Departement de \Mathematiques et d'Informatique},
            city={Trois-{\Rivieres}}, 
            postcode={G9A 5H7},
           %state={},
            country={\Quebec, Canada}}

\begin{abstract}
Bonnet has characterized his surfaces by a geometric condition.
What is done here is a characterization of the same surfaces
by two analytic conditions:
(i) the mean curvature $H$ of a surface in $\mathbb{R}^3$
should admit a reduction to an ordinary differential equation;
(ii) this latter equation should possess the Painlev\'e property.
\end{abstract}

%%Graphical abstract
%\begin{graphicalabstract}
%\includegraphics{grabs}
%\end{graphicalabstract}

%%Research highlights
%  https://www.elsevier.com/researcher/author/tools-and-resources/highlights

%Fading of a script alone does not foster domain-general strategy knowledge

%Performance of the strategy declines during the fading of a script

%Monitoring by a peer keeps performance of the strategy up during script fading

%Performance of a strategy after fading fosters domain-general strategy knowledge

%Fading and monitoring by a peer combined foster domain-general strategy knowledge

%\begin{highlights}
% A purely analytic derivation of Bonnet surfaces
%\item Characterization of Bonnet surfaces by two natural analytic conditions
%\item Research highlight 2
%\end{highlights}

\begin{keyword}
%% keywords here, in the form: keyword \sep keyword
Bonnet surfaces,
Gauss-Codazzi equations,
sixth equation of Painlev\'e.

%% PACS codes here, in the form: \PACS code \sep code

%\noindent \textit{PACS}
\PACS
02.30.Hq, % ordinary differential equations 
02.30.Jr, % Differential equations, partial
02.30.+g  % Mathematical methods in physics
\medskip

%% MSC codes here, in the form: \MSC code \sep code
%% or \MSC[2008] code \sep code (2000 is the default)

%\noindent \textit{MSC} %AMS MSC 2000 Mathematics Subject Classification: 
 % https://mathscinet.ams.org/msc/pdfs/classifications2020.pdf
\MSC
33E17, % Painlev\'e-type functions
34Mxx, % Differential equations in the complex domain [See also 30Dxx, 32G34] 
35A20, % PDEs. Analytic methods, 
35Q99  % Equations of mathematical physics, none of the above
\medskip
    
\end{keyword}

\end{frontmatter}

%% \linenumbers

\bigskip
% https://www.elsevier.com/journals/physics-letters-a/0375-9601/guide-for-authors

TO ELSEVIER.
An error in the style file elsarticle.cls prevents your macro \verb+affiliation+ to accept accented characters
in words like Universit\'e,
Math\'ematiques, Qu\'ebec \dots PLEASE CORRECT. 
In order to overcome this error, these words have been erroneously written without their accents.
See our macros \verb+\Universite+ etc.

%\vfill\eject
\tableofcontents
%% main text

% ============================================================================
\section{Introduction}
\label{section-Bonnet-surfaces}

Surfaces immersed in 
the 
% Euclidean space $\mathbb{R}^3$ 
three-dimensional
simply connected manifold $\mathbb{M}^3(c)$ with constant sectional curvature $\kappa c^2$, 
where $\mathbb{M}^3(c)$ is 
the sphere           $\mathbb{S}^3(c^2)$ when $\kappa= 1$, 
the Euclidean  space $\mathbb{R}^3     $ when $\kappa= 0$ and
the hyperbolic space $\mathbb{H}^3(c^2)$ when $\kappa=-1$, 
are described by a set of partial differential equations (PDE)
known as the Gauss-Codazzi equations\footnote{%
The additional parameter $c$ arises when one replaces
the Euclidean space $\mathbb{R}^3$ 
by the three-dimensional 
manifold $\mathbb{M}^3(c)$ described in the text.
}
\begin{eqnarray}
& &
\left\lbrace 
\begin{array}{ll}
\displaystyle{
\metric_{z \barz} + \frac{1}{2} (H^2-c^2) e^\metric -2 \mod{Q}^2 e^{-\metric}=0 \hbox{ (Gauss)},
}\\ \displaystyle{
Q_\barz-\frac{1}{2} H_z     e^\metric =0,
%}\\ \displaystyle{
\barQ_z-\frac{1}{2} H_\barz e^\metric =0 \hbox{ (Codazzi)},
}
\end{array}
\right.
\label{eqGaussCodazziR3c} 
\end{eqnarray} 
in which 
$z,\barz$ are conformal coordinates,
$e^\metric \D z \D \barz$ the induced metric on the surface,
$H$ the mean curvature ($\metric$ and $H$ being both real),
and $Q$ the complex coefficient of the Hopf differential.

By determining all surfaces in $\mathbb{R}^3$ applicable on a given surface
while conserving the two principal radii of curvature,
the geometer Pierre-Ossian Bonnet \cite{Bonnet1867} isolated a remarkable class of analytic surfaces, 
now called Bonnet surfaces, characterized by
the third order ordinary differential equation (ODE) \cite[Eq.~(52) p.~84]{Bonnet1867}
\cite[Eq.~(B16)]{Conte-Bonnet-JMP}
\begin{eqnarray}
& & {\hskip -16.0truemm}
(\log h')'' + 2 h'
-2\left(\frac{4 \cz}{\sinh 4 \cz (\xi-\xi_0)}\right)^2
  \left(\frac{h'+h^2-c^2}{h'}\right)=0,
	\ '=\frac{\D}{\D \xi},
\label{eq-Bonnet-type4-ODEh}
\end{eqnarray}
in which the mean curvature $H$ is a function of one variable only, denoted $\xi$,
$H(z,\barz)=h(\xi)$, with $\xi=\Re(z)$,
and $\cz^2$ is an arbitrary real constant.
Their metric $\metric$ and the coefficient $Q$ are given by
(see details in \cite{Conte-Lax-PVI-CRAS,Conte-Bonnet-JMP})
% Careful! No x here, only z and \barz
\begin{eqnarray}
& & {\hskip -26.0truemm}
\left\lbrace 
\begin{array}{ll}
\displaystyle{
    Q=2 \cz \coth 2 \cz(    z-    z_0) - 2 \cz  \coth 4 \cz \Re(    z-    z_0)
=\frac{\sinh 2 \cz (\barz-\barz_0)}
      {\sinh 2 \cz (    z-    z_0)}\frac{2 \cz}{\sinh 4 \cz \Re(    z-    z_0)},\
}\\ \displaystyle{
\barQ=2 \cz \coth 2 \cz(\barz-\barz_0) - 2 \cz  \coth 4 \cz \Re(    z-    z_0)
=\frac{\sinh 2 \cz (    z-    z_0)}
      {\sinh 2 \cz (\barz-\barz_0)}\frac{2 \cz}{\sinh 4 \cz \Re(    z-    z_0)},
}\\ \displaystyle{
\modQ^2=\left(\frac{2 \cz}{\sinh 4 \cz \Re(    z-    z_0)}\right)^2,
}\\ \displaystyle{
e^\metric=4 \modQ^2 \frac{\D \Re(z)}{\D H}.
}
\end{array}
\right.
\label{eq-Bonnet-uHQ}
%\label{eq-Bonnet-type4-Q}
\end{eqnarray} 

The general solution \cite{BE1998} of the ODE (\ref{eq-Bonnet-type4-ODEh}),
which admits the first integral 
(\cite[p.~48]{H1897} for $c=0$, \cite{BE2000} for arbitrary $c$)
\begin{eqnarray}
& & {\hskip -15.0 truemm}
K= 
\left(\frac{h''}{h'} + 8 \cz \coth 4 \cz (\xi-\xi_0)\right)^2
\nonumber \\ & & {\hskip -15.0 truemm} \phantom{1234}
+8 \left[
             \left(\frac{4 \cz}{\sinh 4 \cz (\xi-\xi_0)}\right)^2 \frac{h^2 -c^2}{h'}
             + h'
             + 8 \cz \coth 4 \cz (\xi-\xi_0) h\right],
\label{eqHazzi-ODE2h} 
\end{eqnarray}
is the logarithmic derivative of the % one-zero 
tau-function of
either (when $\cz\not=0$) a sixth Painlev\'e equation $\PVI$ of codimension two,  
or     (when $\cz \to 0$) a fifth Painlev\'e equation $\PV $ of codimension three. 
%equivalent to           a codimension-? $\PIII$.

$\PVI$ was only to be discovered  in 1905 \cite{FuchsP6} from purely mathematical considerations,
but in fact it already occured in a natural problem of geometry thirty-eight years earlier,
this is why Bonnet surfaces are so important.

What we perform here a straightforward derivation of
the third order ODE	(\ref{eq-Bonnet-type4-ODEh}),
which is based purely on analysis, and not on geometry like that of Bonnet.

% ============================================================================
\section{The method}
\label{section-The-method}

The mean curvature is an important geometric quantity for a surface in a three-dimensional manifold, therefore
let us make the first 
assumption that this mean curvature $H$ 
depends on one variable denoted $\xi=\varphi(z,\barz)$ instead of two ($z$ and $\barz$),
\begin{eqnarray}
& & H(z,\barz)=h(\xi),\qquad \xi=\varphi(z,\barz),
\label{eqAssumption-H-reduced}
\end{eqnarray}
and let us then determine the function $\varphi(z,\barz)$
and the ODE for $h(\xi)$. 

Let us eliminate $\metric$ and $H$ between the three equations
(\ref{eqGaussCodazziR3c})${}_1$,
(\ref{eqGaussCodazziR3c})${}_2$,
(\ref{eqAssumption-H-reduced})
in order to obtain an equation for $h(\xi)$ whose coefficients 
depend only on $Q,\barQ,\varphi$ and their derivatives.
The field $\metric$ is first determined from the first Codazzi equation (\ref{eqGaussCodazziR3c})${}_2$, 
then substituted in the two remaining Gauss-Codazzi equations (\ref{eqGaussCodazziR3c}). 
After replacement of the field $H$ defined by (\ref{eqAssumption-H-reduced}),
the second Codazzi equation becomes
\begin{eqnarray}
& & Q_z \varphi_z - Q_\barz \varphi_\barz=0,
\end{eqnarray}
and the Gauss equation is now
\begin{eqnarray}
& & {\hskip -6.0truemm}
\left\lbrace 
\begin{array}{ll}
\displaystyle{
(\log h')'' 
+ F_1 h'
+ F_2 \frac{h''}{h'}
+ F_3 \frac{h^2-c^2}{h'}
+ F_4=0,
}\\ \displaystyle{ 
F_1=\frac{Q \barQ}{Q_\barz \varphi_\barz}\ccomma
%}\\ \displaystyle{
F_2=\frac{\varphi_{z \barz}}{\varphi_z \varphi_\barz}\ccomma
%}\\ \displaystyle{
F_3=\frac{Q_\barz}{\varphi_z^2 \varphi_\barz}\ccomma
}\\ \displaystyle{
}\\ \displaystyle{
F_4=
 \frac{Q_{\barz}Q_{\barz\barz}Q_{z\barz}}{\varphi_z   \varphi_\barz}
-\frac{Q_{z\barz\barz} }           {Q_\barz\varphi_z   \varphi_\barz}
-\frac{\varphi_{z\barz} \varphi_{z z}}   {\varphi_z^3 \varphi_\barz}
+\frac{\varphi_{z z \barz}}              {\varphi_z^2 \varphi_\barz}\cdot
}
\end{array}
\right.
\label{eqhnotyetODE}
\end{eqnarray} 
According to our assumption (\ref{eqAssumption-H-reduced}),
equation (\ref{eqhnotyetODE})${}_1$
should be an ODE for $h(\xi)$,
i.e.~its coefficients should only depend on $\xi$.
Let us now make a second assumption, namely that this ODE
should possess the Painlev\'e property \cite{CMBook2},
i.e.~its general solution should have no movable critical singularities.

This double requirement will determine all coefficients $F_j$ (except $F_1$, which will remain arbitrary) in terms of $\xi$.
Once the $F_j(\xi)$'s are known,
the last four equations in (\ref{eqhnotyetODE})
will determine the functions $Q(z,\barz), \barQ(z,\barz), \varphi(z,\barz)$. 

\vfill\eject

% ============================================================================
\subsection{Painlev\'e test} 
\label{sectionPtest}

For the technical vocabulary, we refer the reader to Ref.~\cite[p.~160]{Cargese1996Conte} \cite[Chap.~2]{CMBook2}.

A necessary condition for the absence of movable critical singularities
is that (\ref{eqhnotyetODE}) pass the Painlev\'e test.
Since $h^{(3)}(\xi)$ in (\ref{eqhnotyetODE}) is singular at $h'=0$,
the Painlev\'e test mainly consists
in checking the existence of all Laurent series of
both $h(\xi)$ and $1/h'(\xi)$ 
able to locally represent the general solution.

Near a movable singularity $\chi(\xi)=0$,
since $Q$ is nonzero,
the field $h$ has only one family of movable simple poles,
\begin{eqnarray}
& & h \sim \frac{2}{F_1}\chi^{-1}, \chi \to 0, \chi'=1, 
\end{eqnarray}
and the Fuchs indices  
of this family 
are $-1,1,2$.
The existence of a Laurent series of $h$ near this singularity
requires that two conditions must be obeyed at the positive Fuchs indices,
\begin{eqnarray}
& & 
%Q1=F2(xi)+'(F1,1)(xi)/F1(xi);
%Q2=F1(xi)*F2(xi)**2  + 2* F3(xi) -F1(xi)* F4(xi) + F1(xi)* '(F2 , 1)(xi);
%let rule rQ1 post F2=(function(xi) -deriv(log(F1(xi)),xi)) inactive;
%let rule rQ1a '(F1,1)=(function(xi) -F1(xi)*F2(xi)) inactive;
%let rule rQ2 post F3=(function(xi) F1(xi)*(-F2(xi)**2+F4(xi)-deriv(F2(xi),xi))/2) inactive;
Q_1 \equiv F_1 F_2 + F_1'=0, \qquad
Q_2 \equiv F_1 F_2^2 + 2 F_3 - F_1 F_4 + F_1 F_2'=0. %RC
\label{eqcondQ1Q2}
\end{eqnarray}
They determine $F_2$ and $F_3$ and restrict the ODE (\ref{eqhnotyetODE})${}_1$ to
%
%let rule rF4toG4 post F4=(function(xi) G4(xi)*(F1(xi)/2)**2-deriv(log(F1(xi)),xi,2)+deriv(log(F1(xi)),xi)**2) inactive;
%
%GUESS=deriv(log(deriv(h(xi),xi)),xi,2)+'(h,1)(xi)*F1(xi)
% +((G4(xi)*(F1(xi)/2)**2-deriv(log(F1(xi)),xi,2)+deriv(log(F1(xi)),xi)**2)*'(h,1)(xi)
%   -deriv(log(F1(xi)),xi)*'(h,2)(xi)+(F1(xi)/2)**3*G4(xi)*(h(xi)**2-c**2))/'(h,1)(xi);
%write 0==let rF4toG4 in Bonnetlike1/'(h,1)(xi)**2-GUESS;
\begin{eqnarray}
& & 
\begin{array}{ll}
\displaystyle{
 (\log h')'' + F_1 h' + \frac{F_4 h'+ (F_1/2)^3 G_4 (h^2-c^2)-(\log F_1)' h''}{h'}=0,
}
\end{array}
\label{eqhODEF4}
\end{eqnarray} 
with the notation
\begin{eqnarray}
& & G_4=(2/F_1)^2 \left[F_4 + (\log F_1)'' - {(\log F_1)'}^2\right].
\label{eqrF4toG4}
\end{eqnarray} 

Next, let us build the ODE for $1/h'=g(\xi)$.
After the substitution $h'=1/g(\xi)$ in (\ref{eqhODEF4}) and its derivative,
the elimination of $h(\xi)$ generates a third order ODE for $g(\xi)$,
whose degree in the third derivative is two. Its Gauss decomposition then reads,
%\begin{comment}
%square=(2/F1(xi))**3/G4(xi)*
% (deriv(log(g(xi)),xi,2)/g(xi)-deriv(log(F1(xi)),xi)*deriv(g(xi),xi)/g(xi)**2
%  -F4(xi)/g(xi)-F1(xi)/g(xi)**2+c**2*(F1(xi)/2)**3*G4(xi));
%deriv(square,xi)**2-4*square/g(xi)**2;
%answer*coeff(ODEg,'(g,3)(xi),2)-ODEg*coeff(answer,'(g,3)(xi),2);
%write 0==let rF4toG4 in answer;
%\end{comment}
\begin{eqnarray}
& &
(g S')^2 - 4 S=0,
 S=\frac{8}{F_1^3 G_4} \left[\frac{(\log g)''}{g}- \frac{F_1' g'}{g^3} - \frac{F_4}{g}
- \frac{F_1}{g^2} \right] +c^2=h^2.
\label{eqg}
\end{eqnarray} 
It has one family of movable double poles with an arbitrary leading coefficient,
\begin{eqnarray}
& & g \sim g_0(\xi) \chi^{-2}, \chi \to 0, \chi'=1.
\end{eqnarray}
Its Fuchs indices are $-1,0,2$,
and the index $2$ generates a necessary condition on $G_4, F_1$
for the Laurent series to exist,
\begin{eqnarray}
& & Q_2 \equiv 
\left[\left(\frac{2 G_4'}{F_1 G_4} \right)^2 + 2 G_4 \right]'=0,
\label{eqQ2g}
\end{eqnarray} 
thus defining the value of $G_4$
\begin{eqnarray}
%write let rdcoth in let post G4=(function(xi) 2*k**2*(1-coth(k*FXI(xi))**2)) in first;
%write 0==let '(FXI,1)(xi)=F1(xi)/2 in answer;
& & G_4= 2 (4 \cz)^2 \left[1-\coth^2\left(4 \cz \int (F_1(\xi)/2) \D \xi\right)\right], \cz \hbox{ arbitrary constant}.
\label{eq-value-G4}
\end{eqnarray} 
After substitution of this value (\ref{eq-value-G4}) in (\ref{eqhODEF4}),
the change of independent variable $\xi \to \int (F_1(\xi)/2) \D \xi$
performs the identification of the ODE (\ref{eqhODEF4}) to the ODE (\ref{eqHazzi-ODE2h}) 
for the mean curvature of the Bonnet surfaces.

Therefore, the unique ODE resulting from the two above mentioned assumptions
is the ODE characterizing the Bonnet surfaces.

\vfill\eject

% ============================================================================
\subsection{Determination of the variables $\varphi, Q, \barQ, \metric$} 

Given the values $F_j(\xi)$ obtained in section \ref{sectionPtest},
let us now determine 
the reduced variable $\xi=\varphi(z,\barz)$,
the Hopf coefficients $Q(z,\barz), \barQ(z,\barz)$
and
the metric function $\metric(z,\barz)$
by integrating the last four equations in (\ref{eqhnotyetODE})
and the remaining Gauss-Codazzi equations (\ref{eqGaussCodazziR3c}).

Of course, by the results of Bonnet,
one knows that the expressions (\ref{eq-Bonnet-uHQ}) together with $\varphi=(z+\barz)/2$
are indeed a solution of (\ref{eqhnotyetODE})${}_{2:5}$,
therefore the question is to prove the unicity of this solution.

Without loss of generality, let us first 
make the change of independent variable $\xi \to \int (F_1(\xi)/2) \D \xi$,
i.e.~consider (\ref{eqhnotyetODE})${}_{2:5}$ with the values
given by the conditions (\ref{eqcondQ1Q2}), (\ref{eqrF4toG4}),
\begin{eqnarray}
& & 
F_1=2,
F_2=0,
F_3=F_4=G_4= 2 (4 \cz)^2 \left[1-\coth^2(4 \cz \xi)\right],
\xi=\varphi(z,\barz).
\label{eq-Fj-values}
\end{eqnarray} 

The equation involving $F_2$ in (\ref{eqhnotyetODE})
first yields the function $\varphi$ as a sum of two arbitrary functions of one variable,
\begin{eqnarray}
& & \xi=\varphi=\frac{\Fa(z)+\Fb(\barz)}{2}\cdot
\end{eqnarray} 
The remaining equations (\ref{eqhnotyetODE})${}_{2:5}$ are then equivalent to
\begin{eqnarray}
& & {\hskip -6.0truemm}
\left\lbrace 
\begin{array}{ll}
%    \displaystyle{Q_\barz= - \frac{1}{8} {\Fa'}^2 {\Fb'}   F_4(\xi),
%}\\ \displaystyle{\barQ_z= - \frac{1}{8} {\Fa'} {  \Fb'}^2 F_4(\xi),
%}\\ \displaystyle{Q \barQ= - \frac{1}{8} {\Fa'}^2 {\Fb'}^2 F_4(\xi).
    \displaystyle{Q_\barz= - \frac{1}{8} {\Fa'}^2 {\Fb'}   2 (4 \cz)^2 \left[1-\coth^2\left(4 \cz \frac{\Fa(z)+\Fb(\barz)}{2}\right)\right],
}\\ \displaystyle{\barQ_z= - \frac{1}{8} {\Fa'} {  \Fb'}^2 2 (4 \cz)^2 \left[1-\coth^2\left(4 \cz \frac{\Fa(z)+\Fb(\barz)}{2}\right)\right],
}\\ \displaystyle{Q \barQ= - \frac{1}{8} {\Fa'}^2 {\Fb'}^2 2 (4 \cz)^2 \left[1-\coth^2\left(4 \cz \frac{\Fa(z)+\Fb(\barz)}{2}\right)\right].
}
\end{array}
\right.
\label{eq-sol-remaining}
\end{eqnarray} 
Two quadratures provide $Q$ and $\barQ$ up to an additive arbitrary function
of respectively $z$ and $\barz$,
\begin{eqnarray}
& & {\hskip -6.0truemm}
\left\lbrace 
\begin{array}{ll}
\displaystyle{
    Q=- 2 {\Fa'(    z)}^2 \cz  \coth \left(4 \cz \frac{\Fa(z)+\Fb(\barz)}{2}\right) + f_1(z),
}\\ \displaystyle{
\barQ=- 2 {\Fb'(\barz)}^2 \cz  \coth \left(4 \cz \frac{\Fa(z)+\Fb(\barz)}{2}\right) + f_2(\barz),.
}
\end{array}
\right.
\end{eqnarray} 
The existence of an addition formula for $\coth$ suggests the natural change
of these two arbitrary functions
by
\begin{eqnarray}
& & f_1(    z) =2 {\Fa'(    z)}^2 \cz \coth(2 \cz \tilde f_1(    z)),
    f_2(\barz) =2 {\Fb'(\barz)}^2 \cz \coth(2 \cz \tilde f_2(\barz)).
\end{eqnarray} 
The constraint that the product of these $Q$ and $\barQ$ be equal to (\ref{eq-sol-remaining})${}_3$
is then easily solved as
\begin{eqnarray}
& & \tilde f_1(z)=\Fa(z), \tilde f_2(\barz)=\Fb(\barz).
\end{eqnarray} 
This provides values of $Q$, $\barQ$, $\mod{Q}^2$ identical to the Bonnet values (\ref{eq-Bonnet-uHQ})
upon replacement of, respectively, $\Fa(z)$ by $z-z_0$ and $\Fb(\barz)$ by $\barz-\barz_0$,  
i.e.~by a conformal transformation.
Finally, the value of the metric function $\metric$ follows directly from the Gauss-Codazzi equations.
% ===================================================================
\section{Conclusion}

This direct analytic method could be applied to other problems of differential geometry 
in which one variable is clearly privileged.

\section*{Acknowledgements}

Both authors thank the
Unit\'e Mixte Internationale UMI 3457 du Centre de Recherches
Math\'ematiques de l’Universit\'e de Montr\'eal 
for its financial support. 
This work was partially accomplished during the stay of AMG at the Centre Borelli, 
\'Ecole normale sup\'erieure de Paris-Saclay.
AMG is supported by the Conseil de Recherches en Sciences Naturelles et en G\'enie du Canada (CRSNG).

\vfill\eject

%% If you have bibdatabase file and want bibtex to generate the
%% bibitems, please use
%%
%%  \bibliographystyle{elsarticle-num} 
%%  \bibliography{<your bibdatabase>}

%% else use the following coding to input the bibitems directly in the
%% TeX file.

\vfill\eject
	
\end{document}